\documentclass[11pt]{article}
\PassOptionsToPackage{obeyspaces}{url}
\usepackage{amsfonts, amsmath, amssymb}

\usepackage[hidelinks]{hyperref}
\usepackage{bookmark}
\bookmarksetup{open,numbered,depth=5} % paragraphs go into outlines

\usepackage{amsthm}

\usepackage{tikz}

\usepackage{etoolbox} 
\patchcmd{\thebibliography}{\leftmargin\labelwidth}{\leftmargin\labelwidth\addtolength\itemsep{-0.1\baselineskip}}{}{}

\usepackage{mathtools}

\newcommand*\samethanks[1][\value{footnote}]{\footnotemark[#1]}
\author{Boris Bukh\thanks{Department of Mathematical Sciences, Carnegie Mellon University, Pittsburgh, PA 15213, USA\@. Supported in part by U.S.\ taxpayers through NSF CAREER grant DMS-1555149 and NSF grant DMS-2154063. Email: {\tt bbukh@math.cmu.edu}, {\tt tchao2@andrew.cmu.edu}}
\and
Ting-Wei Chao\samethanks}

\title{Digital almost nets}
\date{}

\usepackage[nameinlink]{cleveref}

\newtheorem{theorem}{Theorem}
\newtheorem{lemma}[theorem]{Lemma}

\newtheorem{definition}[theorem]{Definition}

\newtheorem{claim}{Claim}
\crefname{claim}{Claim}{Claims}

\newcommand*{\eqdef}{\stackrel{\mbox{\normalfont\tiny def}}{=}}   % definition by equality                                      
\newcommand*{\veps}{\varepsilon}                                 % Nice-looking epsilon
\newcommand*{\R}{\mathbb{R}}                                     % Real numbers
\newcommand*{\Z}{\mathbb{Z}}                                     % Integers
\newcommand*{\F}{\mathbb{F}_q}                                   % Finite field of size 2
\newcommand*{\LD}[1]{I_{#1}}                                     % Low-degree polynomials  (Ting-Wei's notation is better)
\DeclarePairedDelimiter\abs{\lvert}{\rvert}                     % Absolute values, cardinality
\DeclareMathOperator{\vol}{vol}                                  % volume (of a box)
\DeclareMathOperator{\rank}{rank}                                % matrix rank
\newcommand*{\T}{\mathcal{T}}                                    % Type of a canonical box
\newcommand*{\A}{\mathcal{A}}                                    % approximation to A (first term of the AP)
\newcommand*{\D}{\mathcal{D}}                                    % approximation to D (step of the AP)
\newcommand*{\Y}{\mathcal{Y}}                                    % approximation to the AP itself
\newcommand*{\bY}{\overline{\mathcal{Y}}}                        % the initial part of \Y
\newcommand*{\cL}{\mathcal{L}}                                   % set of all ways a small box sits inside a canonical box
\newcommand*{\cB}{\mathcal{B}}                                   % box of volume 2^{k-n} (lower bound)

\begin{document}
\maketitle
\begin{abstract}
  Digital nets (in base $2$) are the subsets of $[0,1]^d$ that contain exactly the expected number of points in every not-too-small
  dyadic box. We construct finite sets, which we call ``almost nets'', such that every such dyadic box contains \emph{almost} the
  expected number of points from the set, but whose size is exponentially smaller than the one of nets.
  We also establish a lower bound on the size of such almost nets.
\end{abstract}
\section{Introduction}
We call a subinterval of $[0,1]$ \emph{basic (in base $q$)} if it is of the form $\bigl[\frac{a}{q^k},\frac{a+1}{q^k}\bigr)$, 
for nonnegative integers $a$ and~$k$.
  A \emph{basic box} is a product of basic intervals, i.e.,
  a set of the form $\prod_{i=1}^d\bigl[\frac{a_i}{q^{k_i}},\frac{a_i+1}{q^{k_i}}\bigr)$. 
  If $q=2$, a basic interval is called a dyadic interval, and a basic box is called a dyadic box.

  We say that a set $P\subset [0,1]^d$ is a \emph{$(m,\veps)$-almost net in base $q$} if it is of size $\abs{P}=q^nm$ for some natural
  number $n$ and
\[
  (1-\veps)m\leq \abs{\beta \cap P}\leq (1+\veps)m
\]
for every basic box $\beta$ of volume $\vol(\beta)=q^{-n}$.

In this paper, we are interested in constructions where the parameters $m$ and $\veps$ are independent of~$n$. In contrast, 
since the family of all axis-parallel boxes has finite VC-dimension, one can construct $(m,\veps)$-almost nets
with $m$ linear in $n$ by sampling the points of $P$ at random, see \cite{hs_relative,cm_easy}. The dependence on $n$ is unavoidable
for points sampled at random.

The case $\veps=0$ of the above definition has been well studied. If $\veps=0$, almost nets are known as \emph{digital $(t,m,s)$-nets} or simply \emph{$(t,m,s)$-nets} (in base $q$) 
in the literature\footnote{The parameters $t,m,s$ in the definition of $(t,m,s)$-nets have different meaning than in the present paper. They
correspond to $\log_q m$, $\log_q(mq^n)$ and $d$ respectively in our notation.}.
The adjective `digital' is due to the fact that basic intervals comprise of numbers with specified initial digits in base~$q$.
They are used extensively in discrepancy theory and numerical integration algorithms, and are subject to numerous works, including a book
devoted exclusively to them \cite{dick_pill_book}.
It is known from \cite[Theorem~3]{XN} that, for each $d$, there exist arbitrarily large $(m,0)$-almost nets
with $m\leq q^{5d}$, if $q$ is a prime power. On the other hand, $m$ must grow exponentially with $d$ for large enough nets \cite{martin_visentin} (see also \cite{schurer} for asymptotic analysis of
the bound in \cite{martin_visentin}).

In contrast to these results, for $\veps>0$, we construct $(m,\veps)$-almost nets with $m$ being only polynomial in~$d$. 

\begin{theorem}\label{mainthm}
For any prime $q$, any $d\geq 2$, and any positive integers
$m,n$ satisfying $m\geq 400d\log (dq)$, there exists a set $P\subset [0,1]^d$ of size $mq^n$ such that, for any basic box $\beta$ of volume $q^{-n}$, 
\begin{equation}\label{eq:main}
  \left(1-10\sqrt{\frac{d\log (dq)}{m}}\right)m\leq \abs{\beta\cap P}\leq\left(1+10\sqrt{\frac{d\log (dq)}{m}}\right)m.
\end{equation}
In particular, for every $0<\veps<1/2$ and every $d\geq 2$, there exist arbitrarily large $(m,\veps)$-almost nets in base $q$ with $m\leq 100\veps^{-2}d\log(dq)$.

Furthermore, the set $P$ satisfying \eqref{eq:main} can be chosen to be an $\bigl(M,0)$-net in base $q$ with $M\leq\nobreak d^{4d}q^{6d}m$.
\end{theorem}
This result has an application in geometric Ramsey theory: A \emph{convex hole} in a finite set $S\subset \R^d$
is a subset $H\subset S$ in convex position and whose convex hull contains no other point of~$S$.
An old problem of Valtr \cite{valtr} asks for the largest $h(d)$ such that every sufficiently large $S\subset \R^d$
in general position contains a hole of size~$h(d)$. Using \Cref{mainthm} one can show that $h(d)\leq 4^{d+o(d)}$, which is
an improvement over the bound of $h(d)\leq 2^{7d}$ that can be obtained from~$(t,m,s)$-nets. The 
details of both bounds are in \cite{bukh_chao_holzman}.

The construction behind \Cref{mainthm} is a minor modification on the construction in \cite{bukh_chao}. Whereas
the construction in \cite{bukh_chao} uses primes in $\Z$, this construction uses irreducible polynomials
in $\F[x]$. The reason for this change is to make the denominators be powers of the same prime~$q$.
Furthermore, because the addition in $\F[x]$ satisfies the ultrametric inequality (with respect to the degree), and because we do not need to
worry about boxes that are not basic, several details in the new construction are simpler. As such, we do not make
any claims about the novelty. Our purpose in writing the present note is to record the details
of the construction for its application to convex holes. We also hope that almost nets will
find applications in many other areas that currently use the conventional nets.\medskip

We do not know when the bound in \Cref{mainthm} is sharp. The following is the best lower bound we were able to prove.
Its dependence on $\veps$ is close to optimal, as long as $\veps$ is not too small, but the dependence on $d$ is poor.
In the special case $\veps=0$, we recover the lower bound $t= \Omega(s)$ in $(t,m,s)$-nets via a proof
different than those in \cite{martin_visentin,schurer} (keeping in mind that that $t$ and $s$ correspond
to $\log_q m$ and $d$ respectively).

\begin{theorem}\label{lowerbound}
Assume that there exists an $(m,\varepsilon)$-net $P\subseteq [0,1]^d$ in base $q$, then the following holds. 

If $\varepsilon\geq 1/2\sqrt{d}$, then
\[m= \Omega\bigl(\frac{\log d}{q^2\varepsilon^2\log(1/\varepsilon)}\bigr).\]

If $1/2\sqrt{d}\geq \varepsilon\geq e^{-d/8}$, then 
\[m= \Omega(q^{-2k-2}\varepsilon^{-2}),\]
where $k=\frac{2\log(1/\varepsilon)}{\log d-\log\log(1/\varepsilon)}$.

In particular, if $\varepsilon=\omega(d^{-t})$ for some constant $t$, then we have $m= \Omega_{q,t}(1/\varepsilon^2)$.

If $\varepsilon=o(e^{-cd})$ for some constant $c$ such that $0<c<\min(1/8,1/q^2)$, then we get an exponential lower bound 
$m=\Omega\bigl(q^{-2}e^{c'd}\bigr)$, where $c'=2c(1-2\log q/\log(1/c))$.
\end{theorem}

\paragraph{Open problem.} It would be interesting to prove a result similar to \Cref{mainthm} which applies to all boxes, not only to basic boxes. It is possible to construct a set $P$ for which
$\abs{\beta\cap P}$ is lower-bounded by $(1-o(1))m$ for all boxes $\beta$ by taking
a union of translates of the set $P$ from \Cref{mainthm} in a manner similar to
that in the second part of the proof of \cite[Theorem~2]{bukh_chao}. However, we
have been unable to control $\abs{\beta\cap P}$ from above.

\paragraph{Acknowledgment.} We are thankful to Ron Holzman for useful discussions, and to two anonymous referees for their help in improving
this paper.

\section{Proof of \texorpdfstring{\Cref{mainthm}}{Theorem 1}}
We denote by $\F$ the finite field consisting of elements $0,1,2,\dotsc,q-1$ equipped with the usual
mod-$q$ arithmetic. Let $t\eqdef \lceil 2\log_q d+2\rceil$. Since the number of irreducible polynomials of degree $t$ in $\F[x]$ is
\[ 
\frac{1}{t}(\sum_{i|t}\mu(i)q^{t/i})\geq \frac{1}{t}(q^t-q^{t/2+1})\geq d,
\]
we may pick $d$ distinct irreducible polynomials $p_1,\dotsc,p_d$ of degree $t$ in $\F[x]$. We fix some such choice
of polynomials for the duration of the proof. We associate each of these $d$ polynomials to the respective coordinate direction.
We will be interested in \emph{canonical boxes}, which are the boxes of the form
\[
  B=\prod_{i=1}^d \left[\frac{a_i}{q^{k_it}},\frac{a_i+1}{q^{k_it}}\right).
\]
for some nonnegative integers $k_i$ and $0\leq a_i<q^{k_it},i=1,2,\dotsc,d$.

We say that a polynomial $f\in \Z[x]$ is a \emph{basic polynomial} if $\deg f<t$ and all of its coefficients are in~$\{0,1,\dotsc,q-1\}$.

For an irreducible polynomial $p\in \F[x]$ of degree $t$ and a polynomial $f\in\F[x]$, we define the \emph{base-$p$ expansion of $f$} to be
$f=\nobreak f_0+f_1p+\dotsb+f_{\ell}p^{\ell}$, where each $f_i$ is a basic polynomial.
Put $r_p(f)\eqdef (f_0(1/q)+f_1(1/q)q^{-t}+\dotsb+f_{\ell}(1/q)q^{-\ell t})/q$, where we view the basic polynomials
$f_0,f_1,\dotsc,f_{\ell}$ as polynomial functions on $\R$. 
In other words, if $f_i=\sum_{j<t} c_{i,j}x^j$ with $c_{i,j}\in\{0,1,\dotsc,q-1\}$,
then, denoting by $C_i$ the concatenation $c_{i,0}c_{i,1}\ldots c_{i,t-1}$, the base-$q$ expansion of the real
number $r_p(f)$ is
\[
  r_p(f)=0.C_0C_1\ldots C_{\ell}.
\]
Note that $r_p(f)\in [0,1)$. 
Define the function $r\colon \F[x]\to [0,1]^d$ by $r(f)\eqdef \bigl(r_{p_1}(f),\dotsc,r_{p_d}(f)\bigr)$.

Recall that our aim is to construct a set $P\subset [0,1]^d$ whose intersection with any basic box of volume $q^{-n}$
has almost the expected number of points.
\begin{definition}
We say that a box $\beta$ is \emph{good} if $\beta$ is a basic box of volume $\vol(\beta)=q^{-n}$. 
Let $B$ be the smallest canonical box containing~$\beta$. We call $(B,\beta)$ a \emph{good pair}.
\end{definition}
Note that if $(B,\beta)$ is a good pair, then $\vol(B)\leq q^{-n+dt-1}$. Indeed, every basic interval
is contained in an interval of the form $[a/q^{kt},(a+1)/q^{kt})$ that is at most
$q^{t-1}$ times larger, and \linebreak therefore $\vol(B)\leq q^{-n}(q^{t-1})^d\leq q^{-n+dt-1}$.\medskip

Suppose $B$ is a canonical box. Write it as $B=\prod_i\left[a_i/q^{k_it},(a_i+1)/q^{k_it}\right)$, and consider $r^{-1}(B)$. The set $r^{-1}(B)\subseteq \F[x]$ consists
of all solutions to the system
\begin{align*}
f&\equiv a_1'\pmod{p_1^{k_1}},\\
f&\equiv a_2'\pmod{p_2^{k_2}},\\
% From https://tex.stackexchange.com/questions/7650/centering-vdots-in-a-system-of-many-equations
 &\setbox0\hbox{$\equiv$}\mathrel{\makebox[\wd0]{\vdots}}\\
f&\equiv a_d'\pmod{p_d^{k_d}},
\end{align*}
where % $a_i'\eqdef \sum_{j=0}^{k_i-1}f_{i,j}p_i^j$
$a_i'=f_{i,0}+f_{i,1}p_i+\dotsb+f_{i,k_i-1} p_i^{k_i-1}$
and $f_{i,0},f_{i,1},\ldots,f_{i,k_i-1}$ are the unique basic polynomials satisfying $a_i/q^{k_it}=(f_{i,0}(1/q)+f_{i,1}(1/q)q^{-t}+\ldots+f_{i,k_i-1}(1/q)q^{-(k_i-1)t})/q$. 

By the Chinese Remainder theorem, the set $r^{-1}(B)$ is of the form $A(B)+D(B)\F[x]$ where $D(B)\eqdef p_1^{k_1}p_2^{k_2}\dotsb p_d^{k_d}$ and $A(B)$ is the unique element in $r^{-1}(B)$ of degree less than $t(k_1+\ldots+k_d)$. Note that $\deg D(B)=t(k_1+\ldots+k_d)=-\log_q (\vol(B))$.

Given a good pair $(B,\beta)$, define \[L_B(\beta)\eqdef\{ g \in\F[x] : r\bigl(A(B)+gD(B)\bigr) \in \beta\}.\]

\begin{claim}\label{claim:L}
The set $\cL\eqdef \lbrace L_B(\beta):(B,\beta)\mbox{ is a good pair}\rbrace$ is of size at most $q^{4dt}$.
\end{claim}

\begin{proof}
Let $(B,\beta)$ be a good pair. Write $B$ and $\beta$ in the form
\[B=\prod_{i=1}^d \left[\frac{a_i}{q^{k_it}},\frac{a_i+1}{q^{k_it}}\right),\qquad \beta=\prod_{i=1}^d \left[\frac{a_i}{q^{k_it}}+\frac{b_i}{q^{(k_i+1)t}},\frac{a_i}{q^{k_it}}+\frac{c_i}{q^{(k_i+1)t}}\right).\]

The condition $r\bigl(A(B)+gD(B)\bigr)\in \beta$ is equivalent to
\begin{align*}
A(B)+gD(B)&\in a_1'+p_1^{k_1}J_1\pmod{p_1^{k_1+1}},\\
A(B)+gD(B)&\in a_2'+p_2^{k_2}J_2\pmod{p_2^{k_2+1}},\\
&\setbox0\hbox{$\in$}\mathrel{\makebox[\wd0]{\vdots}}\\
A(B)+gD(B)&\in a_d'+p_d^{k_d}J_d\pmod{p_d^{k_d+1}},
\end{align*}
where the sets $J_i$ consist of the basic polynomials $f$ such that $f(1/q)q^{t-1}\in [b_i,c_i)$.

On the other hand, 
\begin{align*}
A(B)+gD(B)&\equiv a_1'+(\alpha_1+g\delta_1)p_1^{k_1}\pmod{p_1^{k_1+1}},\\
A(B)+gD(B)&\equiv a_2'+(\alpha_2+g\delta_2)p_2^{k_2}\pmod{p_2^{k_2+1}},\\
&\setbox0\hbox{$\equiv$}\mathrel{\makebox[\wd0]{\vdots}}\\
A(B)+gD(B)&\equiv a_d'+(\alpha_d+g\delta_d)p_d^{k_d}\pmod{p_d^{k_d+1}}
\end{align*}
for some $\alpha_i,\delta_i\in\F[x]/(p_i),i=1,2,\dots,d$. Since $\dim_{\F} \F[x]/(p_i)=\deg p_i=t$, there are at most $q^{2dt}$ different choices for $(\alpha_i,\delta_i)_{i=1}^d$.
Also, there are at most $q^{2dt}$ different choices for $(b_i,c_i)_{i=1}^d$ satisfying $0\leq b_i<c_i\leq q^t$. Since $L_B(\beta)$ is determined by $(\alpha_i,\delta_i,b_i,c_i)_{i=1}^d$, the claim is true.
\end{proof}

To each canonical box $B$ of volume between $q^{-n}$ and $q^{-n+dt-1}$ inclusive we assign
a \emph{type}, so that boxes of the same type behave similarly.
Formally, let $\A(B)$ be the polynomial obtained from the polynomial $A(B)$ by setting the coefficients of $1,x,x^2,\dotsc,x^{n-dt-1}$ to zero. Similarly, let $\D(B)$ be the polynomial obtained from $D(B)$
by setting the coefficients of $1,x,x^2,\dotsc,x^{n-3dt}$ to zero. The type of $B$ is then the pair $\T(B)\eqdef \bigl(\A(B),\D(B)\bigr)$. 

Note that, from $q^{-n}\leq \vol(B)\leq q^{-n+dt-1}$ and $\deg D(B)=-\log_q(\vol(B))$ it follows that
\begin{equation}\label{eq:dbound}
  n-dt+1\leq\deg \D(B)\leq n.
\end{equation}

\begin{claim}\label{claim:type}
The number of types is at most $q^{4dt}$.
\end{claim}
\begin{proof}
Since $\deg A(B)<\deg D(B)\leq n$, only the $dt$ (resp.~$3dt$) leading coefficients of $\A(B)$ (resp.~$\D(B)$) may be non-zero. Hence, the number of types is at most $q^{dt}\times q^{3dt}=q^{4dt}$.
\end{proof}

For a type $\T=(\A,\D)$, let $\Y(\T)\eqdef \lbrace \A+g\D : g\in\F[x]\rbrace$. Note that if $\T=\T(B)$, then $\Y(\T)$ is an approximation to $r^{-1}(B)$.
That is to say, the respective elements of $\Y(T)$ and of $r^{-1}(B)$ differ only in low-degree coefficients.

Let $\LD{k}$ denote polynomials of degree less than $k$ in $\F[x]$.
Our construction will be a union of sets of the form $h+\LD{n-dt}$ where $\deg h\leq n+dt$.\smallskip

We first prove that there is no difference in how the sets $\Y(\T)$ and $r^{-1}(B)$ intersect $h+\LD{n-dt}$.
\begin{claim}\label{claim:shrink}
Suppose $\T(B)=\bigl(\A(B),\D(B)\bigr)$.
Then for any polynomial $h\in \LD{n+dt}$ and any polynomial $g$,
$\A(B)+g\D(B)\in h+\LD{n-dt}$ if and only if $A(B)+gD(B)\in h+\LD{n-dt}$.
\end{claim}

\begin{proof}
If $\A(B)+g\D(B)\in h+\LD{n-dt}$, then $\deg (\A(B)+g\D(B))< n+dt$. Since $\deg \A(B)< n$ and $\deg \D(B)\geq n-dt$, it follows that $\deg g< 2dt$. From the definition of $\A(B)$ and $\D(B)$, the coefficients of $x^{n-dt},x^{n-dt+1},\ldots$ in $\A(B)+g\D(B)$ are the same as the respective coefficients in $A(B)+gD(B)$. The opposite direction is similar.
\end{proof}

For a type $\T$ and $L\in\cL$ that satisfy $\T=\T(B)$ and $L=L_B(\beta)$ for some good pair $(B,\beta)$, define
\[\Y_{\T}(L)\eqdef\lbrace \A+g\D:g\in L\rbrace.\]
With this definition,
$\Y_{\T}(L)$ 
is the approximation to $r^{-1}(\beta)$ induced by the approximation $\Y(\T)$
to~$r^{-1}(B)$.

\begin{claim}\label{clm:ex}
The set $\bY_{\T}(L)\eqdef\Y_{\T}(L)\cap \LD{n+dt}$ is of size exactly $q^{dt}$.
\end{claim}

\begin{proof}
Let $(B,\beta)$ be a good pair such that $\T=\T(B)$ and $L=L_B(\beta)$. From the previous claim, we know that the size of $\bY_{\T}(L)$ is the same as the size of $r^{-1}(\beta)\cap \LD{n+dt}$. By the Chinese remainder theorem, each of the canonical boxes of volume $q^{-(\lfloor n/t\rfloor+d)t}$ contains equally many points from $r(\LD{n+dt})$. Since $n\leq (\lfloor n/t\rfloor+d)t$, the number of points in $\beta\cap r(\LD{n+dt})$ is equal to $q^{n+dt}\vol(\beta)=q^{dt}$.
\end{proof}

\begin{claim}
  Let $h$ be chosen uniformly from $\LD{n+dt}$. Then
  $\abs{\bY_{\T}(L)\cap (h+\LD{n-dt})}$ is $1$ with probability $q^{-dt}$ and is $0$ otherwise.
\end{claim}
\begin{proof}
  Let $u\in \bY_{\T}(L)$ be arbitrary. Clearly $\Pr[u\in h+\LD{n-dt}]=q^{-2dt}$. The events of the form $u\in h+\LD{n-dt}$ are mutually disjoint as $u$ ranges over $\bY_{\T}(L)$. Indeed,
  suppose $\T=(\A,\D)$ and $u,u'\in \bY_{\T}(L)$ are such that $u,u'\in h+\LD{n-dt}$ for some $h\in \LD{n+dt}$. We may write $u=\A+g\D$ and $u'=\A+g'\D$. Then $u-u'=(g-g')\D\in \LD{n-dt}$. Since $\deg \D(B)\geq n-dt$, this implies that $g=g'$ and hence $u=u'$.

In the combination with \Cref{clm:ex}, this implies that 
\[\Pr\bigl[\,|\bY_{\T}(L)\cap (h+\LD{n-dt})|=1 \bigr]=q^{-2dt}q^{dt}=q^{-dt}.\qedhere\]
\end{proof}

Sample $q^{dt}m$ elements uniformly at random from $\LD{n+dt}$, independently from one another.
Let $H$ be the resulting multiset, and consider the multiset $H+\LD{n-dt}\eqdef \{h+f : h\in H,f\in \LD{n-dt}\}$. For a type $\T$ and $L\in\cL$ that satisfy $\T=\T(B)$ and $L=L_B(\beta)$ for some good pair $(B,\beta)$, define the random variable $N_{\T,L}\eqdef |\bY_{\T}(L)\cap (H+\LD{n-dt})|$.
This random variable is distributed according to the binomial distribution $\operatorname{Binom}(q^{dt}m,q^{-dt})$.

Let $\veps=\sqrt{33dt\log q/m}$. Note that $\veps<\sqrt{33d(2\log d+3\log q)/m}<10\sqrt{d\log (dq)/m}$, and in particular
$\veps<1/2$. Hence, $\veps^2/2-\veps^3/2\geq \veps^2/4$. By the tail bounds for the binomial distribution \cite[Theorems A.1.11 and A.1.13]{alonspencer} we obtain
\[\Pr\bigl[\,N_{\T,L}-m>\varepsilon m \bigr]<e^{-(\varepsilon^2/2-\varepsilon^3/2)m}<q^{-8dt}/2,\]
\[\Pr\bigl[\,N_{\T,L}-m<-\varepsilon m\bigr]<e^{-\varepsilon^2m/2}<q^{-8dt}/2.\]

From \Cref{claim:L,claim:type} and the union bound it then follows that there exists a choice of $H$
such that $N_{\T,L}$ is bounded between $(1-\veps)m$ and $(1+\veps)m$ whenever $\T=\T(B)$, $L=L_B(\beta)$ and $(B,\beta)$ is a good pair. By \Cref{claim:shrink}, this implies that the number of points in any good box $\beta$ of volume $q^{-n}$, the size $\beta\cap r(H+\LD{n-dt})$ is bounded between $(1-\veps)m$ and $(1+\veps)m$.

Hence the multiset $r(H+\LD{n-dt})$ in $[0,1)^d$ is of size exactly $mq^n$ and satisfies \eqref{eq:main}. Since the $r$\nobreakdash-image of
every set of the form $h+\LD{n-dt}$, for $h\in \F[x]$, is a $(q^{dt},0)$-net, it follows that $r(H+\LD{n-dt})$ is a $(M,0)$-net with $M=q^{2dt}m\leq d^{4d}q^{6d}m$.

To obtain a set satisfying the same conclusion, we may perturb the points of $r(H+\LD{n-dt})$ slightly to ensure distinctness.

\section{Proof of \texorpdfstring{\Cref{lowerbound}}{Theorem 2}}
We shall derive \Cref{lowerbound} from the following lemma. 
\begin{lemma}\label{lemma:lower}
For any positive integers $n,d, q$ and positive real numbers $m,\varepsilon$ with $n\geq d\geq 2$, $q\geq 2$, and $\veps<1/4$, if there exists an $(m,\veps)$-almost net $P\subseteq [0,1]^d$ in base $q$ of size $q^nm$, then 
\[m= \Omega\bigl(\frac{\log(\binom{d}{k})}{q^{2k}\varepsilon^2\log(1/\varepsilon)}\bigr),\]
for any integer $k$ such that $1\leq k\leq d/2$ and 
\begin{equation}\label{eq:epscond}
2\varepsilon\geq \binom{d}{k}^{-1/2}
\end{equation}
holds.
\end{lemma}
\begin{proof}
Let $\cB$ be the box $[0,1/q^{n-2k})\times [0,1)^{d-1}$. For any point $v=(v_1,\ldots,v_d)\in \cB$, 
write its coordinates in base $q$ as $v_\ell=(0.v_{\ell,1}v_{\ell,2}\ldots)_q$. Noting that
the first $n-2k$ base-$q$ digits of $v_1$ are zero, we let $X_1(v)$ be the first non-trivial
digit of $v_1$, i.e., $X_1(v)\eqdef v_{1,n-2k+1}$. Similarly, let $X_\ell(v)\eqdef v_{\ell,1}$ for $\ell\geq 2$.

The proof idea is to use almost independence of functions $X_1,\dotsc,X_d$ for a randomly chosen
point of $\cB$. However, we do not directly appeal to the known bound on the size of probability
spaces supporting almost independent random variables (see e.g. \cite{alon,aakmrx}) because
those bounds are formulated for $\{0,1\}$-valued random variables, whereas $X_1,\dotsc,X_d$
take $q$ distinct values. \medskip

Let $S\eqdef P\cap \cB$, and $t\eqdef \abs{S}$. Since $P$ is an $(m,\veps)$-almost net, it follows that $t$ is between $q^{2k}(1-\varepsilon)m$ and $q^{2k}(1+\varepsilon)m$. Assume $v^1,\ldots,v^t$ are all the points in $S$.

For $x\in \R$, let $e_q(x)\eqdef \exp(2\pi ix/q)$ where $i=\sqrt{-1}$.
Let $U$ be a $\binom{d}{k}$-by-$t$ matrix, where the rows are indexed by $\binom{[d]}{k}$ and the columns are indexed by $[t]$. The general entry of $U$ is
\[
  U_{J,\ell}\eqdef e_q\bigl(\sum_{j\in J}X_j(v^\ell)\bigr).
\]
Also, define $A\eqdef \tfrac{1}{t}UU^*$.

\begin{claim}
The diagonal terms in $A$ are all $1$. The off-diagonal terms are, in absolute value, bounded above by $2\varepsilon$.
\end{claim}
\begin{proof}
The general term of $A$ is given by 
\[
  A_{J_1,J_2}=\tfrac{1}{t}\sum_{\ell=1}^te_q\Bigl(\sum_{j\in J_1}X_j(v^\ell)-\sum_{j\in J_2}X_j(v^\ell)\Bigr).
\]
If $J_1=J_2$, this is clearly $1$.

Suppose $J_1\neq J_2$. Note that, for any choice of $\alpha=(\alpha_j)_{j\in J_1\Delta J_2}$ with $\alpha_j\in \{0,1,\dotsc,q-1\}$, the set 
\[
  %B_{\alpha}\eqdef 
  \{v\in \cB : X_j(v)=\alpha_{j}\text { for }j\in J_1\Delta J_2\}
\]
is a basic box of volume $q^{-n+2k-|J_1\Delta J_2|}$. Thus, for any $\tau\in \{0,1,\dotsc,q-1\}$, the region
\[B_\tau \eqdef \Bigl\{v\in\cB:\sum_{j\in J_1}X_j(v)-\sum_{j\in J_2}X_j(v)\equiv\tau\pmod q\Bigr\}\]
can be partitioned into $q^{|J_1\Delta J_2|-1}$ many basic boxes of volume $q^{-n+2k-|J_1\Delta J_2|}$ each. Since we have $\break-n+2k-|J_1\Delta J_2|\geq -n$, it follows that the number of $\ell$ such that $v^\ell\in B_\tau$ is bounded between $q^{2k-1}(1-\varepsilon)m$ and $q^{2k-1}(1+\varepsilon)m$. Thus, 
\begin{align*}
|A_{J_1,J_2}|=&\frac{1}{t}\abs*{\sum_{\tau=1}^q \abs{B_\tau\cap S}e_q(\tau)}\\
\leq& \frac{1}{t}\abs*{\sum_{\tau=1}^q q^{2k-1}m e_q(\tau)}+\frac{1}{t}\sum_{\tau=1}^q |\varepsilon q^{2k-1}me_q(\tau)|\\
=&\frac{\varepsilon q^{2k}m}{t}\leq 2\varepsilon.\qedhere
\end{align*}
\end{proof}

We apply \cite[Theorem~2.1]{alon} to the matrix $(A+\bar{A})/2$. We obtain that, if 
$
\binom{d}{k}^{-1/2}\leq 2\varepsilon < 1/2,
$ 
then $2q^{2k}(1+\varepsilon)m\geq 2\rank(A)\geq \rank\bigl((A+\bar{A})/2\bigl)= \Omega(\frac{\log(\binom{d}{k})}{\varepsilon^2\log(1/\varepsilon)})$. Therefore,
\[m= \Omega\bigl(\frac{\log(\binom{d}{k})}{q^{2k}\varepsilon^2\log(1/\varepsilon)}\bigr).\qedhere\]
\end{proof}

The right hand side of lemma \ref{lemma:lower} is a decreasing function of $k$ for $k\in [1,d/2]$. Therefore, we shall pick $k$ as small as possible. If $\varepsilon\geq 1/2\sqrt{d}$, then we may set $k=1$ and get
\[m= \Omega\bigl(\frac{\log d}{q^2\varepsilon^2\log(1/\varepsilon)}\bigr).\]
If $1/2\sqrt{d}\geq \varepsilon\geq e^{-d/8}$, then we may set $k=\frac{2\log(1/\varepsilon)}{\log d-\log\log(1/\varepsilon)}$. From the assumption on~$\varepsilon$, we have $k\leq \log(1/\varepsilon)$. Therefore, 
\begin{align*}
\binom{d}{k}\geq& (d/k)^k\\
\geq& \exp\bigl(\frac{2\log(1/\varepsilon)}{\log d-\log\log(1/\varepsilon)}(\log d-\log k)\bigr)\\
\geq& \exp\bigl(\frac{2\log(1/\varepsilon)}{\log d-\log\log(1/\varepsilon)}(\log d-\log\log(1/\varepsilon))\bigr)\\
\geq& \frac{1}{\varepsilon^2},
\end{align*}
and so \eqref{eq:epscond} holds.
Hence, we may apply \Cref{lemma:lower} with $\lceil k\rceil$ in place of $k$ and obtain
\[m= \Omega(q^{-2k-2}\varepsilon^{-2}).\]
In particular, if $\varepsilon=\omega(d^{-t})$ for some constant $t$, then $k$ is also a constant, and so $m= \Omega_{q,t}(1/\varepsilon^2)$ in this case.

If $\varepsilon=o(e^{-cd})$ for some constant $c$ such that $0<c<\min(1/8,1/q^2)$, then the $(m,\varepsilon)$-net is also an $(m,e^{-cd})$-net, when $d$ is large enough. We may apply the result above with $e^{-cd}$ in place of $\varepsilon$. In this case, the calculations above yield $k=2cd/\log (1/c)$, and we get 
$m= \Omega\bigl(q^{-2}e^{c'd}\bigr)$ where $c'=2c\bigl(1-2\log q/\log(1/c)\bigr)$.

\bibliographystyle{plain}
\bibliography{AlmostNet}
\end{document}